\documentclass[12pt]{article}
\usepackage{amsmath}
\usepackage{amsfonts}
\usepackage{amssymb}
\usepackage{latexsym}
\usepackage{graphicx}
\usepackage{enumerate}
\usepackage[shortlabels]{enumitem}
\usepackage{mathtools}
\usepackage{hyperref}

\usepackage{amsmath, amssymb,latexsym}
\usepackage{color}
\usepackage{xcolor}
\usepackage{fullpage}

\def\C{\mathbb{C}}
\def\F{\mathbb{F}}
\def\la{\lambda}

\def\T{\mathsf{T}}

\def\e{\mathbf{e}}

\DeclareMathOperator{\diag}{diag}

\DeclareMathOperator{\rank}{rank}
\DeclareMathOperator{\im}{Im}

\newtheorem{theorem}{Theorem}[section]

\newtheorem{lemma}[theorem]{Lemma}
\newtheorem{definition}[theorem]{Definition}

\title{Computing a compact local Smith McMillan form}
\author{
Vanni Noferini\thanks{Aalto University, Department of Mathematics and Systems Analysis, P.O. Box 11100, FI-00076, Aalto, Finland. Supported by an Academy of Finland grant (Suomen Akatemian p\"{a}\"{a}t\"{o}s 331240).}
\and 
Paul Van Dooren\thanks{Universit\'e catholique de Louvain, Department of Mathematical Engineering, Av. Lemaitre 4, B-1348 Louvain-la-Neuve, Belgium. Supported by an Aalto Science Institute Visitor Programme.}
}
\begin{document}
\maketitle
\begin{abstract} 
We define a compact local Smith-McMillan form of a rational matrix $R(\la)$ as the diagonal matrix whose diagonal elements are the nonzero entries of a local Smith-McMillan form of $R(\la)$. We show that a recursive rank search procedure, applied to a block-Toeplitz matrix built on the Laurent expansion of
$R(\la)$ around an arbitrary complex point $\la_0$, allows us to compute a compact local Smith-McMillan form of that rational matrix $R(\la)$ at the point $\la_0$,
provided we keep track of the transformation matrices used in the rank search. It also allows us to recover the root polynomials 
of a polynomial matrix and root vectors of a rational matrix, at an expansion point $\la_0$. Numerical tests illustrate the promising performance of the resulting algorithm.
\end{abstract}



\section{Introduction}

Finding poles and zeros of a rational matrix $R(\la)\in \C(\la)^{m\times n}$ with coefficients in the field of complex numbers $\C$ is one of the basic problems in linear system theory. Such a rational matrix describes the input/output behaviour of a general system of differential or difference equations \cite{kailath,Rosenbrock}. Its poles correspond to the natural frequencies of the dynamical system, while its zeros correspond to the frequencies that are blocked by the system \cite{MFK}. When $R(\la)$ does not have full row or column rank over the field $\C(\lambda)$ of rational functions, the rational matrix $R(\la)$ has also a non trivial left (respectively, right) null space \cite{For75}, which yields additional information on the initial conditions and degrees of freedom related to the response of the system \cite{For75,kailath,MFK,Rosenbrock}.
A rational matrix can have multiple poles and zeros in a point $\la_0\in \C$ and even coalescent poles and zeros. The finer structure of the response of the dynamical system at the frequency $\la_0$ that is such a pole/zero  is then described by the local Smith-McMillan form of the rational matrix. The latter is a diagonal matrix that associates with each pole/zero a number of  structural indices that reflect the structure of the response of the system at that frequency.

Finding poles and zeros and their structural indices is therefore an important problem in the analysis of a dynamical system. However, the response of the system also depends on particular directions. These are
input vectors that either excite one of the system's poles or are blocked by one of the system's zeros, of a particular multiplicity given by each structural index.
In the latter case, these {\em root vectors} \cite{DopN21,GLR82,N11,N23,NofV22,NofV22a} arise naturally when one wants to describe the solution set of particular matrix equations involving rational matrices \cite{kailath,kar94} or appear as expansion vectors in tangential interpolation problems of high order \cite{GVV04}. Root vectors can be viewed as a generalization of an eigenvector for a first order system of differential or difference equations modeled by the eigenvalue problem $\la x- Ax=0$. The structural indices at that zero are then linked to its Jordan structure. 
In this paper we show how to compute such a local Smith-McMillan form at a pole/zero $\la_0$ of $R(\la)$ by applying the Toeplitz rank search algorithm \cite{vdv} to a block Toeplitz matrix built on the Laurent expansion of $R(\la)$ around $\la_0\in \C$. We also link this method to the calculation of root polynomials of a polynomial matrix, or root vectors of a rational matrix, at an expansion point $\la_0$ \cite{DopN21,N11,NofV22,NofV22a}.

The paper is organized as follows. In Section \ref{Sec:background}, we recall the basic definitions and background material for the rest of the paper and we introduce the compact local Smith-McMillan form. In Section \ref{Sec:Compact} we recall the rank properties of triangular Toeplitz matrices defined from the Laurent expansion at a point $\la_0\in \C$ and show their relation to the computation of the structural indices of a rational matrix. Section \ref{Sec:Compact} contains the main new results of the paper: it shows that the Toeplitz rank search also constructs a local Smith form of a polynomial matrix and a local Smith-McMillan form of a rational matrix.  In Section \ref{Sec:numerics}, we show some numerical experiments indicating that the accuracy of the algorithm is very satisfactory. Finally, we give some concluding remarks in Section \ref{Sec:conclusion}.

\section{Background and definitions} \label{Sec:background}

\subsection{The local Smith-McMillan form}
We denote the field of rational functions with complex coefficients by
$\C(\lambda)$ and the ring of polynomials with complex coefficients by $\C[\lambda]$. 
The structure at a finite point $\la_0\in \C$ which is a pole or a zero of an $m \times n$ rational matrix
$R(\la)\in \C(\la)^{m\times n}$ is defined via its 
local Smith-McMillan form at the point $\la_0\in \C$ \cite{MCM}:
\begin{equation} \label{MCM}
 M(\la)R(\la)N(\la) := 
\left[ \begin{array}{ccc|c} (\la-\la_0)^{\sigma_1} & & 0 & \\
 & \ddots & & \\ 0 & &  (\la-\la_0)^{\sigma_r}& \\ \hline
 & & & 0_{m-r,n-r} \end{array}\right] ,
\end{equation}
where $M(\la)$ and $N(\la)$ are rational, well defined and such that
$M(\la_0),N(\la_0)$ are defined and invertible (that is, $\la_0$ is neither a zero nor a pole of $M(\la),N(\la)$), $r$ is the normal rank of $R(\la)$, i.e.\ the rank of $R(\la)$ over the field $\C(\la)$, and the integers $\sigma_i$ are called {\it structural indices} of $R(\la)$ at the point
$\la_0$, and are ordered non-decreasingly, i.e., $\sigma_1 \leq \ldots \leq \sigma_r.$ 
The negative indices refer to poles of the transfer function and the positive indices refer to zeros of the transfer function. An index $\sigma_j=0$ is not associated with any dynamical behaviour and corresponds to neither a pole nor a zero. The standard Smith-McMillan form has the same structure, but the transformation matrices $M(\la)$ and $N(\la)$ are then unimodular and the diagonal elements are then the so-called invariant factors $e_i(\la)/f_i(\la)\in \C(\la)$ where the polynomials $e_i(\la)\in \C[\la]$ and $f_i(\la)\in \C[\la]$ are monic and satisfy the divisibility chains
$$ e_1(\la) | e_2(\la) | \ldots| e_r(\la), \quad \mathrm{and}
\quad f_r(\la) |\ldots | f_2(\la) | f_1(\la) $$
(see e.g. \cite{MCM}). 
We point out that the local form can be derived from the standard Smith-McMillan form via a one-sided extraction
of the factors $(\la-\la_0)^{\sigma_i}$ from the invariant factors  $e_i(\la)/f_i(\la)$, which implies that we can choose one of the two
matrices $M(\la)$ and $N(\la)$ in the local form \eqref{MCM} to be polynomial and unimodular. The classical computation of the standard Smith-McMillan decomposition is based on the Euclidean algorithm and Gaussian elimination over the ring of polynomials, which precludes numerical pivoting techniques and is therefore numerically unreliable \cite{vdv}.
For this reason it has been suggested as a better alternative to compute it via linearizations \cite{NofV22,NofV22a}. In this paper we show that the local decomposition
can also be obtained from the Laurent expansion around the point $\la_0$.

\subsection{Null spaces and their minimal indices}

When defining the structure of a general $m\times n$ rational matrix, one typically includes 
the structure of its right null space $\ker R(\la)$ and left null space $\ker R(\la)^T$,
which are rational vector spaces over the field $\C(\la)$. 
Their characterization is based on  particular polynomial bases, for which we need the following definition.
\begin{definition} 
A matrix polynomial $N(\la) \in \C[\la]^{n\times p}$ of normal rank $p$
is called a minimal polynomial basis if the sum of the degrees of its columns, 
called the order of the basis, is minimal among all polynomial bases of the range of $N(\la)$, i.e., the vector space of all $\C(\la)$-linear combinations of the columns of $N(\la)$. Its ordered column degrees are called the minimal indices of the basis.
\end{definition}
It is known \cite{For75} that the ordered list of indices is independent of 
the choice of minimal basis of the space. One can define the right null space $\ker R(\la)$ and left null space 
$\ker R(\la)^T$ of an $m\times n$  rational matrix $R(\la)$ of normal rank $r$ as 
the vector spaces of rational vectors $x(\la)$ and $y(\la)$ annihilated by $R(\la)$ on the respective sides:
$$  \ker R(\la):=\{ x(\la) \; | \; R(\la)x(\la)=0\}, \qquad    \ker R(\la)^T :=\{ y(\la) \;| \; y^\T(\la) R(\la)=0\}.
$$
Then, the minimal indices of any minimal polynomial basis for these spaces, are called the right and left  
minimal indices of $R(\la)$. The dimensions of $\ker R(\la)$ and $\ker R(\la)^T$ are, respectively, $n-r$ and 
$m-r$ and the right and left minimal indices are denoted by
$$\{\epsilon_1,\ldots,\epsilon_{n-r}\}, \quad \{\eta_1,\ldots,\eta_{m-r}\}.
$$
It is known \cite{For75} that for any minimal basis $N(\la)$, the constant matrix 
$N(\la_0)$ has full column rank for all $\la_0\in\C$ and the highest column degree matrix
of $N(\la)$ also has full column rank.

\subsection{The compact local Smith-McMillan form and the compact local Smith form}

The local Smith-McMillan form not only contains information on the structural indices $\sigma_i, i=1,\ldots,r$, but the
invertible matrices $M(\la)$ and $N(\la)$ in \eqref{MCM} also contain bases for the left and right null spaces of the rational matrix $R(\la)$. It follows indeed from \eqref{MCM} that 
$$  \ker R(\la) = \im \left( N(\la)\left[\begin{array}{c} 0 \\ I_{n-r}\end{array}\right] \right), \quad 
 \ker [R(\la)]^\mathsf{T} = \im  \left( [M(\la)]^\mathsf{T}\left[\begin{array}{c} 0 \\ I_{m-r}\end{array}\right] \right).
$$ 
These block columns are invertible bases, in the sense of \cite{N23}, of the modules, $\ker R(\la)\cap \C[\la]$ and $\ker [R(\la)]^T\cap \C[\la]$), respectively, because they have full column rank for all finite $\la=\la_0 \in \C$ (see \cite[Corollary 4.2]{N23} for an analogous argument in the case of analytic matrix functions), but they are not necessarily minimal polynomial bases.

\smallskip

The following more compact equation discards the information of the null spaces and focuses only on the structural indices at $\lambda_0$.
We call the diagonal matrix containing the nonzero local invariant factor a {\em compact local Smith-McMillan form}, in analogy to the compact SVD of a matrix $A$, which also discards the singular values and vectors related to the left and right null space of the constant matrix $A$. We say that a rational matrix $R(\la)$ is left (resp. right) invertible at $\la_0$ if $\la_0$ is not a pole and the constant matrix $R(\la_0)$ is left (resp. right) invertible.
\begin{theorem} \label{compactMCM}
Every rational matrix $R(\la)\in \C^{m\times n}(\la)$ of normal rank $r$ satisfies the following equations revealing its compact local Smith-McMillan form:
\begin{equation} \label{MCMr}
R(\la)N_r(\la) =
\hat M_r(\la)\left[ \begin{array}{ccc} (\la-\la_0)^{\sigma_1} & & 0 \\
 & \ddots & \\ 0 & &  (\la-\la_0)^{\sigma_r} \end{array}\right] ,
\end{equation}
where $\sigma_1,\dots,\sigma_r$ are the structural indices of $R(\la)$ at $\la_0$, $N_r(\la)\in \C[\la]^{n\times r}$ has a polynomial left inverse and $\hat M_r(\la)\in \C(\la)^{m\times r}$ is left invertible at $\la_0$, and
\begin{equation} \label{MCMl}
M_\ell(\la)R(\la) =
\left[ \begin{array}{ccc} (\la-\la_0)^{\sigma_1} & & 0 \\
 & \ddots & \\ 0 & &  (\la-\la_0)^{\sigma_r} \end{array}\right] \hat N_\ell(\la)
\end{equation}
where $\sigma_1,\dots,\sigma_r$ are the structural indices of $R(\la)$ at $\la_0$, $M_\ell (\la)\in \C[\la]^{r\times m}$ has a polynomial right inverse and $\hat N_\ell (\la)\in \C(\la)^{r\times n}$ is right invertible at $\la_0$.
\end{theorem}
\begin{proof}
We only give a proof for the first form \eqref{MCMr} since the form \eqref{MCMl} is dual to it. We start from the decomposition
\eqref{MCM} where we can make the choice that $N(\la)$ is unimodular, and $M(\la)$ is rational and such that $M(\la_0)$ is invertible. If we then multiply it on the left with $\hat M(\la):= M^{-1}(\la)$ and define
$$
 N_r(\la) := N(\la)\left[\begin{array}{c} I_r \\ 0\end{array}\right], \quad 
 \hat M_r(\la) = \hat M(\la) \left[\begin{array}{c} I_r \\ 0 \end{array}\right],
$$
then the result follows, since the property of $N_r(\la)$ follows from the unimodularity of $N(\la)$, and the existence and left invertibility of $\hat M_r(\la_0)$ follows from the invertiblity of $M(\la_0)$. \hfill
\end{proof}

If $R(\la)$ happens to be polynomial, then $\hat M_r(\la)$ and $\hat N_\ell(\la)$ can also be chosen to be polynomial, which then yields a compact local Smith form, as stated in the following theorem.
\begin{theorem}
 \label{compactSM}
Every polynomial matrix $P(\la)\in \C^{m\times n}[\la]$ of normal rank $r$ satisfies the following equations revealing its compact local Smith form:
\begin{equation} \label{SMr}
P(\la)N_r(\la) =
\hat M_r(\la)\left[ \begin{array}{ccc} (\la-\la_0)^{\sigma_1} & & 0 \\
 & \ddots & \\ 0 & &  (\la-\la_0)^{\sigma_r} \end{array}\right] ,
\end{equation}
where  $\sigma_1,\dots,\sigma_r$ are the structural indices of $R(\la)$ at $\la_0$, $N_r(\la)\in \C[\la]^{n\times r}$ has a polynomial left inverse and $\hat M_r(\la)\in \C[\la]^{m\times r}$ is left invertible at $\la_0$, and
\begin{equation} \label{SMl}
M_\ell(\la)P(\la) =
\left[ \begin{array}{ccc} (\la-\la_0)^{\sigma_1} & & 0 \\
 & \ddots & \\ 0 & &  (\la-\la_0)^{\sigma_r} \end{array}\right] \hat N_\ell(\la)
\end{equation}
where $\sigma_1,\dots,\sigma_r$ are the structural indices of $R(\la)$ at $\la_0$, $M_\ell (\la)\in \C[\la]^{r\times m}$ has a polynomial right inverse and $\hat N_\ell(\la)\in \C[\la]^{m\times r}$ is right invertible at $\la_0$.
\end{theorem}

\subsection{Connection with root vectors and root polynomials} \label{Subs:root}

The Smith form and the Smith-McMillan form, respectively, are closely related to the concepts of (left and right) root polynomials of a polynomial matrix \cite{DopN21,GLR82,N11,NofV22a} and (left and right) root vectors of a general rational matrix \cite{NofV22}.
The definition for the right vectors requires an invertible basis \cite{N23}, i.e., an arbitrary polynomial basis for the right null space of the matrix that, upon evaluation at $\la=\la_0$, spans the same subspace as a minimal basis. For instance, $N(\la)\begin{bsmallmatrix} 0 \\ I_{n-r}
\end{bsmallmatrix} $, where $N(\la)$ is as in \eqref{MCM}, is such a basis; indeed it is an example of what was called an \emph{invertible basis} of a module in \cite{N23} as even though it may not be a minimal basis, it nevertheless has full rank upon evaluation at any finite point. It is then useful to introduce some notation to denote those column vectors of the matrices $N_r(\la)$ and $\hat M_r(\la)$, in Theorems \ref{compactMCM} and \ref{compactSM},  that correspond to the positive structural indices~:
$$ j:=\min(i \; | \; \sigma_i\ge1), \quad x_i(\la):=N_r(\la)\e_i, \; j \le i \le r, \quad v_i(\la):=\hat M_r(\la)\e_i, \; j \le i \le r.$$
The following properties of the vectors $x_i(\la)$ follow from the compact local Smith form at the zero $\la_0$ of a polynomial matrix $P(\la)$.
The vectors $x_i(\la)$ satisfy the equations 
\begin{equation*} %\label{rootpoly}
P(\la)x_i(\la) = v_i(\la)(\la-\la_0)^{\sigma_i},  \quad  \ker P(\la_0)=\im \left(\left[ \begin{array}{cccc} N(\la_0)\begin{bsmallmatrix} 0 \\ I_{n-r}
\end{bsmallmatrix}  & x_j(\la_0)& \ldots & x_r(\la_0)\end{array}\right]\right)
\end{equation*}
where the matrices 
\begin{equation} \label{rp}
\left[ \begin{array}{cccc} N(\la_0)\begin{bsmallmatrix} 0 \\ I_{n-r}\end{bsmallmatrix}  & x_j(\la_0)& \ldots & x_r(\la_0)\end{array}\right] \quad \mathrm{and} \quad \left[ \begin{array}{ccc}v_j(\la_0)& \ldots & v_r(\la_0)\end{array}\right]
\end{equation}  
have full column rank; in particular the rank of the matrix on the left of \eqref{rp} is equal to $\dim \ker P(\la_0)$. These properties of the vectors $x_i(\la)$ are precisely the defining properties of a \emph{complete set of \it root polynomials} of $P(\la)$ at $\la_0$, and they follow directly from the local compact Smith form
at $\la_0$. Moreover, one  can also show that such vectors are maximal sets of root polynomials \cite{DopN21}.

\medskip

For a rational matrix $R(\la)$, one again looks only at the positive structural indices $\sigma_i\ge 1$ and the same definition holds for  the column vectors $x_i(\la)$ and $v_i(\la)$, and again we have
\begin{equation*} %\label{rootvec}
R(\la)x_i(\la) = v_i(\la)(\la-\la_0)^{\sigma_i},  \quad  \ker R(\la_0)=\im \left( \left[ \begin{array}{cccc} N(\la_0)\begin{bsmallmatrix} 0 \\ I_{n-r}
\end{bsmallmatrix}  & x_j(\la_0)& \ldots & x_r(\la_0)\end{array}\right] \right)
\end{equation*} 
where still the full rank conditions of \eqref{rp} hold. These properties also follow directly from the compact local Smith-McMillan form. Again, one can show that the $x_i(\la)$ are maximal sets of root vectors, see \cite{NofV22}. We point out in particular that $\ker R(\la_0)$ can still be defined even when $\la_0$ is a pole \cite[Definition 3.8]{NofV22}, and it does not contain the directions in which $R(\la)$ tends to infinity when $\la \rightarrow \la_0$.

\medskip

The link with root polynomials and root vectors is one of the main motivations to construct a compact Smith-McMillan form. While in the proof of Theorem \ref{compactMCM} $N_r(\la)$ and $\hat{M}_r(\la)$ were constructed starting from a full local Smith-McMillan form, algorithmically it is more efficient to compute a compact local Smith-McMillan form directly, as opposed to computing the full one and only later discard some columns. This leads to the question of whether the rightmost columns of a matrix $N_r(\la)$ satisfying \eqref{MCMr} are still a complete set of root vectors, i.e., satisfy the rank condition in \eqref{rp}.
The following Lemma implies that they do.

\begin{lemma} \label{Lem:extend}
Let $R(\la)N_r(\la)=\hat M_r(\la)\Lambda(\la)$ with 
$\Lambda(\la):=\diag((\la-\la_0)^{\sigma_1},\ldots,(\la-\la_0)^{\sigma_r})$ be a compact local Smith-McMillan form as defined in Theorem \ref{compactMCM}, and let $\tilde N(\la) \in \C[\la]^{m\times (n-r)}$ be any completion such that $U(\la):=\left[N_r(\la), \; \tilde N(\la) \right]$ is unimodular. Then, there exists a polynomial matrix $Y(\la)$ such that a polynomial basis for $\ker R(\la)$ is given by $B(\la):=\tilde N(\la)-N_r(\la)Y(\la)$. Moreover, $B(\la_0)$ spans the same space as $N(\la_0) \begin{bmatrix}
    0\\I_{n-r}
\end{bmatrix}$ as in \eqref{rp}.
\end{lemma}
\begin{proof}
Without loss of generality, let us suppose $R(\la)\neq 0$. Let $\pi(\la)$ be any nonzero scalar polynomial such that $\pi(\la)R(\la)$ is polynomial; for example, we can take $\pi(\la)=f_1(\lambda)$, the denominator of the $(1,1)$ element in the Smith-McMillan form of $R(\la)$. By assumption, $\pi(\la) R(\la) U(\la)$ is unimodularly equivalent over $\C[\la]$ with  $\begin{bmatrix} \pi(\la) R(\la)N_r(\la) & 0\end{bmatrix}$.  We can then invoke Theorem \ref{ring}, whose statement and proof we postpone to the Appendix, to conclude that $\pi(\la) R(\la) \tilde N(\la) = \pi(\la) R(\la) N_r(\la) Y(\la)$ for some polynomial matrix $Y(\la)$. 
On the other hand, $B(\la) = U(\la) \begin{bmatrix} -Y(\la)\\
I_{n-r} \end{bmatrix}$, and hence $\rank B(\la) = n-r$
and therefore $B(\la)$ is a basis for $\ker R(\la)=\ker \pi(\la)R(\la)$. Finally, to prove the last statement, define $\hat{N}(\la)=N(\la) \begin{bmatrix}
    0\\
    I_{n-r}
\end{bmatrix}$ where $N(\la)$ is the unimodular matrix appearing in \eqref{rp}. Then, both $B(\la)$ and $\hat{N}(\la)$ are invertible polynomial bases \cite{N23} for the same $\C[\la]$-module, i.e., $\ker R(\la) \cap \C[\la]$. This implies that $\hat{N}(\la_0)$ and $B(\la_0)$ are both full rank. We conclude that $B(\la_0)$ and $\hat{N}(\la_0)$ span the same $\C$-vector subspace, i.e., $\ker_{\la_0} R(\la)$ as defined in \cite{DopN21,N11,NofV22}.
\end{proof}

  \smallskip  

  \noindent
In particular, it follows from Lemma \ref{Lem:extend} that the first of the rank conditions
\eqref{rp} hold whenever the vectors $x_i(\la)$ are the rightmost columns of a matrix $N_r(\la)$ satisfying Theorem \ref{compactMCM}. The second rank condition then follows from the first by properties of maximal sets of root vectors \cite{NofV22}. 

\bigskip

To conclude this section, we note that one can give definitions for the left root vectors or root polynomials that are dual to those of the right vectors, and use the left compact local decompositions instead. Details are therefore left out.

\section{Constructing a compact local Smith-McMillan form} \label{Sec:Compact}
In this section, we describe an algorithm to compute the compact local Smith-McMillan form as in Theorem \ref{compactMCM}.
\subsection{Retrieving the structural indices} \label{Sec:StructuralIndices}

We first recall here an important connection between the structural indices $\{\sigma_i,i=1,\ldots,r\}$ of a pole/zero $\la_0\in \C$ of a
general rational matrix $R(\la)$ and its Laurent expansion around that point \cite{vdv}. Let us assume that the pole $\la_0$ is of order $\ell$, and that it is possibly also a zero. Then $R(\la)$ has a Laurent expansion about the point $\la_0$, with leading coefficient $R_{-\ell}~$:
\begin{equation} \label{Laurent}
R(\la) := R_{-\ell}(\la-\la_0)^{-\ell} + R_{-\ell +1}(\la-\la_0)^{-\ell+1} + R_{-\ell+2}
(\la-\la_0)^{-\ell+2} + R_{-\ell+3}(\la-\la_0)^{-\ell+3} + \ldots \end{equation}
The following theorem derives the structural indices at $\la_0$ from the expansion \eqref{Laurent}.
\begin{theorem}[\cite{vdv}] \label{th:rank}
Using the coefficients of the Laurent expansion \eqref{Laurent}, let us define for $k\ge -\ell$, the block Toeplitz matrices 
\begin{equation*} %\label{BlockToeplitz}
T_{\la_0,k}(R) := \left[ \begin{array}{cccc}
R_{-\ell} & R_{-\ell+1} & \ldots & R_k \\
& R_{-\ell} & \ddots & \vdots \\
& & \ddots & R_{-\ell+1} \\
& & & R_{-\ell}
\end{array}\right] \in \C^{m(k+\ell+1)\times n(k+\ell+1)}.
\end{equation*}
Let their ranks and rank increments be denoted by $r_k:=\rank T_{\la_0,k}(R)$, and $\rho_k:= r_k-r_{k-1}$, where we set
$r_k=0$ for $k < -\ell$. Then the number $e_i$ of indices $\sigma_j$ that are equal to $i$, is given by
\begin{equation*} %\label{rankTkpoles}  
 e_i:=\#\{ \sigma_j=i\} = \rho_i-\rho_{i-1}= r_i-2r_{i-1}+r_{i-2}.
\end{equation*}  
Moreover, the surplus ranks $\rho_k$ form a non-decreasing sequence
\[  0 \le \rho_{-\ell} \le \ldots \le \rho_{d'} = r  \] 
and $d':=\max_i(\sigma_i)$ is the smallest index $k$ for which $\rho_k=r$, the normal rank of $R(\la)$.
\end{theorem}
For simplicity, when no ambiguity arises we will denote the Toeplitz matrices $T_{\la_0,k}(R)$ by just $T_k(R)$ or $T_k$.
It follows from the above theorem that one only has to compute the ranks  of the sequence
$\{T_k, -\ell \le k \le d'\}$, and hence, one only needs to know the coefficients $\{R_k, -\ell \le k \le d'\}$
of the expansion. If $d'=\max_i \sigma_i$ is not known in advance, we will see that this index is also discovered by the algorithm provided the normal rank of $R(\la)$ is known. The latter can be estimated, for example, by evaluating the rank of the transfer function in some randomly generated points.

\subsection{Toeplitz rank search} \label{Sec:RankSearch}

In \cite{vdv} a {\em Toeplitz rank search} algorithm was proposed to compute the rank increments of Theorem \ref{th:rank}, while exploiting the block Toeplitz structure of the matrices $T_{k}(R)$. In this paper we slightly modify this algorithm so that it also 
constructs a compact local Smith-McMillan decomposition, by keeping track of the intermediate transformations.

To simplify the derivation, we first consider the case where $R(\la)$ does not have a pole at the finite point $\la_0$ but only a zero. 
Then $R(\la)$ has a Taylor expansion at that point
\[
R(\la) := R_0 + R_1(\la-\la_0) +
R_2(\la-\la_0)^2 + R_3(\la-\la_0)^3 + \ldots \]
and the corresponding Toeplitz matrices then have the leading coefficient $R_0$ on diagonal 
\begin{equation} \label{Toeplitz} T_k := \left[ \begin{array}{cccc}
R_0 & R_1 & \ldots & R_k \\
& R_0 & \ddots & \vdots \\
& & \ddots & R_1 \\
& & & R_0
\end{array}\right], \quad r_k:= \rank T_k \quad \mathrm{for}\: k\ge 0.
\end{equation}
Note that polynomial matrices
are a special case of such rational matrices having no poles at $\la_0$, since all their poles are at infinity.
Moreover, the Toeplitz rank search for the polynomial matrix $P(\la)$ obtained by truncating the Taylor expansion of the rational matrix $R(\la)$ after its 
first $d'+1$ coefficients $\{R_i, 0 \le i \le d'\}$ has the same structural indices, according to Theorem \ref{th:rank}.
We therefore focus first on polynomial matrices.

We recall the algorithm derived in \cite{vdv} for computing the structural indices of a rational matrix in a pole/zero at $\la_0$ from its Laurent expansion. We apply it here to the expansion about $(\la-\la_0)$ of a polynomial matrix
\begin{equation*} %\label{P}
P(\la)=P_0+ P_1(\la-\la_0) + \ldots + P_d(\la-\la_0)^d. 
\end{equation*}
To simplify our notation, we will assume in this section that $\la_0=0$. This does not affect the generality of the results.
This algorithm computes a rank factorization of the Toeplitz matrices $T_k$ given in \eqref{Toeplitz}.
It operates on the stacked array of coefficients $P_i$ using a sequence of invertible transformations, followed by shifts of sub-blocks~:
\begin{equation*} %\label{recur}
\left[ \begin{array}{cc}
 L_d^{(0)} & R_d^{(0)} \\  L_{d-1}^{(0)} & R_{d-1}^{(0)} \\ \vdots & \vdots \\ L_1^{(0)} & R_1^{(0)} \\  L_0^{(0)} & 0
\end{array}\right] := \left[ \begin{array}{cc}
 P_d \\ P_{d-1} \\ \vdots \\ P_1 \\ P_0
\end{array}\right] N_0 ,  \qquad   \left[ \begin{array}{cc}
 P^{(0)}_d  \\ P^{(0)}_{d-1} \\ \vdots \\ P^{(0)}_1 \\ P^{(0)}_0
\end{array}\right]  :=  \left[ \begin{array}{cc}
 L_d^{(0)} & 0 \\  L_{d-1}^{(0)} & R_{d}^{(0)} \\ \vdots & \vdots \\ L_1^{(0)} & R_2^{(0)} \\  L_0^{(0)} & R_1^{(0)} 
\end{array}\right]     
\end{equation*}
where $N_0$ is an invertible column transformation compressing the columns of $P_0$,  $L_0^{(0)}$ has full column rank $r_0=\rho_0$, and nullity $n_0=\nu_0:=n-\rho_0$, and $L_i^{(0)}\in\C^{m\times \rho_0}$ and $R_i^{(0)}\in\C^{m\times \nu_0}$. 

\smallskip

If we apply the same invertible transformation to block columns of the Toeplitz matrix $T_d$ then we obtain 
$$
\left[ \begin{array}{cccccccc}
L_0^{(0)} & 0 & L_1^{(0)} & R_1^{(0)} & \ldots  & L_d^{(0)} & R_d^{(0)} \\
& & L_0^{(0)} & 0  & \ddots & \vdots & \vdots \\
& & & &   \ddots & L_1^{(0)} & R_1^{(0)}  \\
& & & & & L_0^{(0)} & 0 
\end{array}\right] :=
 \left[ \begin{array}{ccccc}
P_0 & P_1 & \ldots  & P_d \\
& P_0 & \ddots & \vdots \\
& &  \ddots & P_1 \\
& & & P_0
\end{array}\right]  \left[ \begin{array}{ccccc}
N_0  \\
& N_0  \\
& &  \ddots \\
& & & N_0
\end{array}\right] ,
$$
and after permuting the $\left[ L_i^{(0)} \; R_i^{(0)}\right] $ pairs to  $\left[R_i^{(0)} \; L_i^{(0)}\right]$ 
pairs, this becomes
\begin{equation} \label{submatrix}
\left[ \begin{array}{cccccc}
0 & \! P^{(0)}_0 & \ldots  & P^{(0)}_{d-1} \! & L_d^{(0)} \\
& &  \ddots & \vdots &\vdots \\
& & & P^{(0)}_0 & L_1^{(0)} \\
& & & & L_0^{(0)}
\end{array}\right]  \! := \! 
 \left[ \begin{array}{ccccccc}
0 & \! L_0^{(0)} \! & R_1^{(0)} & L_1^{(0)} & \ldots  & R_d^{(0)} & \! L_d^{(0)} \\
& & 0 & L_0^{(0)}  & \ddots & \vdots & \vdots \\
& & & &   \ddots & R_1^{(0)} & \! L_1^{(0)}  \\
& & & & & 0 & \! L_0^{(0)} 
\end{array}\right] .
\end{equation}
This shows that $\rank T_d(P) = \rho_0 + \rank T_{d-1}(P^{(0)})$, and is the basis of a recursive computation of the successive ranks of a block Toeplitz matrix $T_d$ and its submatrices. We repeat this
on the polynomial matrix $P^{(0)}(\la)$ and its corresponding Toeplitz matrix $T_{d-1}(P^{(0)})$ which turns out to be a submatrix of the left hand side of \eqref{submatrix}. This induction step is repeated on the subsequent polynomial matrices $P^{(k)}(\la)$, and shows that we finally compress the column space of $T_d$ by induction using an invertible transformation $N$ that is a product of the individual invertible block-diagonal transformations and permutations~:
$$  \left[ \begin{array}{ccccc}
P_0 & P_1 & \ldots  & P_d \\
& P_0 & \ddots & \vdots \\
& &  \ddots & P_1 \\
& & & P_0
\end{array}\right]  =
 \left[ \begin{array}{cccccccc}
0 & L_0^{(d)} & L_1^{(d-1)} & \ldots  & L_d^{(0)} \\
& & L_0^{(d-1)}  &  & \vdots \\
& & &   \ddots & L_1^{(0)}  \\
& & & & L_0^{(0)} 
\end{array}\right]\cdot N, $$
where the suffixes $^{(i)}$ refer to the iteration step $i$ of the Toeplitz rank search.
Here each ``diagonal" block $L_0^{(i)}$ has rank $\rho_i$ and 
$$ 0 \le \rho_0 \le \rho_1 \le \ldots \le \rho_d.
$$
These rank inequalities follow easily from the above algorithmic construction since $L_0^{(i+1)}$ is a column compression of 
the compound matrix $[L_0^{(i)},R_1^{(i)}]$. It also follows from this that 
$$  \rank T_k(P) = \sum_{i=0}^k\rho_i $$
and it was shown in \cite{vdv} that $\rho_i= \# \{ \sigma_j\le i\}$ and $e_i:= \rho_i-\rho_{i-1}= \# \{\sigma_j =  i\}$.

\medskip

In order to link this to a compact local Smith form, we write these operations as a polynomial matrix equation (i.e. where 
$P^{(0)}(\la) $ is polynomial as well)~:
\begin{equation} \label{factorize}
 P(\la)N_0 = P^{(0)}(\la)  \Lambda_0(\la), \quad \mathrm{where} \quad  \Lambda_0(\la):= \left[ \begin{array}{cc}
 I_{\rho_0} & 0 \\ 0 & \la I_{\nu_0}
\end{array}\right].
\end{equation}
We will show that the first $\rho_0$ columns of \eqref{factorize} already match those of the compact local Smith form.
Let us now look at the factorization after the next step, yielding 
\begin{equation} \label{factor2}
 P^{(0)}(\la)N_1 = P^{(1)}(\la) \Lambda_1(\la), \quad \mathrm{where} \quad  \Lambda_1(\la):=\left[ \begin{array}{cc} I_{\rho_1} & 0 \\ 0 & \la I_{\nu_1} \end{array}\right]
\end{equation}
It follows from \eqref{factorize} and \eqref{factor2} that 
$$  P(\la) N_0 \Lambda_0^{-1}(\la)N_1\Lambda_0(\la)=  P^{(1)}(\la) \Lambda_0(\la)\Lambda_1(\la), \;  \mathrm{where} \;  
\Lambda_0(\la)\Lambda_1(\la) =\left[ \begin{array}{ccc} I_{\rho_0}  \\  & \la I_{\rho_1-\rho_0} \\ & & \la^2 I_{\nu_1} \end{array}\right].
$$
This clearly goes in the right direction, provided the matrix $\Lambda_0^{-1}(\la)N_1\Lambda_0(\la)$ is unimodular.
In \cite{vdv} the Toeplitz rank search was implemented with unitary transformations $N_i$ in order to guarantee good numerical stability properties. This allowed to reconstruct the partial multiplicities $\sigma_j$ at the considered root, but 
if one also wants to reconstruct a compact local Smith form, then one also needs to satisfy the different conditions described in Theorems \ref{compactMCM} and \ref{compactSM}. Therefore one needs to constrain the rank search to a special set of transformations, as explained below. The column rank compression          
$$ \left[\begin{array}{cc} L^{(i-1)}_0 & R^{(i-1)}_1 \end{array}\right] N_{i} = \left[\begin{array}{c|cc}  L^{(i-1)}_0 & L^{(i-1)}_{0+} &0 \end{array}\right] =: \left[\begin{array}{c|cc} L^{(i)}_0 & 0 \end{array}\right] 
$$
where $L_0^{(i-1)}$ has full rank $\rho_{i-1}$ and $L^{(i)}_0:=[L^{(i-1)}_{0} \; L^{(i-1)}_{0+}\;]$ has full rank 
$\rho_{i}\ge \rho_{i-1}$, can be implemented as a factored transformation
with a simple inverse
\begin{equation} \label{Mi} N_{i}= \left[\begin{array}{cc} I_{\rho_{i-1}} & Z_i \\ 0 & I_{\nu_{i-1}} \end{array}\right]
\left[\begin{array}{cc} I_{\rho_{i-1}} & 0  \\ 0 & Q_i  \end{array}\right]= 
\left[\begin{array}{cc} I_{\rho_{i-1}} & Z_iQ_i \\ 0 & Q_i  \end{array}\right] , \quad N^{-1}_{i}=
 \left[\begin{array}{cc}  I_{\rho_{i-1}} & -Z_i \\ 0 & Q^*_i  \end{array}\right] ,
\end{equation} 
where $Z_i=-{L_0^{(i-1)}}^\dagger R_1^{(i-1)}$ is the least squares solution of $L_0^{(i-1)}Z_i=-R_1^{(i-1)}$ and $Q_i$ is a unitary transformation compressing the columns of $(I_{\nu_{i-1}}-{L_0^{(i-1)}}{L_0^{(i-1)}}^\dagger)R_1^{(i-1)}$ to
the $(\rho_{i}-\rho_{i-1})$ independent columns of $L_{0+}^{(i)}$ that, by construction, are also independent from those of $L_0^{(i)}$. If the matrix $R_1^{(i-1)}$ lies in the span of $L_0^{(i-1)}$, then $\rho_{i}=\rho_{i-1}$ and the matrix $Q_i=I_{\nu_{i-1}}$.
 Note also that in the special case $i=0$ the equations also hold with $R_1^{(-1)}=P_0$, $\rho_{-1}=0$ and $\nu_{-1}=n$.
 For this particular choice of transformations, we now have the following result.
\begin{theorem} \label{th:fact}
The choice of transformations $N_i$ given in \eqref{Mi} for the Toeplitz rank search algorithm produces the factorization 
\begin{equation} \label{factor}
P(\la) N(\la)= P^{(d')}(\la)\Lambda(\la), \quad \Lambda(\la):=\Lambda_0 \cdots\Lambda_{d'}, \quad N(\la):= N_0\Lambda_0^{-1} N_1\Lambda_1^{-1} \cdots N_{d'}\Lambda_{d'}^{-1}\Lambda(\la)
\end{equation}
where $N(\la)$ is unimodular,  
$$\Lambda(\la)=\diag(I_{e_0}, \la I_{e_1},\ldots,\la^{d'-1} I_{e_{d'-1}},\la^{d'} I_{\nu_{d'}}), $$
$e_i:=\rho_i-\rho_{i-1}$ for $i\ge 0,$ and the constant matrix $P^{(d')}(0)= \left[L_0^{(d')}, \; 0 \;\right]$ has rank $r$ and nullity $n-r$.
\end{theorem}
\begin{proof}
It follows from the recursive rank search algorithm that it stops as soon as
$$ P(\la)N_0\Lambda_0^{-1} N_1\Lambda_1^{-1} \cdots N_{d'}\Lambda_{d'}^{-1}= P^{(d')}(\la).
$$ with  $P^{(d')}(0)= \left[L_0^{(d')}, \; 0 \;\right]$ having rank $r$.
If we multiply both sides with $\Lambda(\la)$, then we also obtain
$$ P(\la) N_0\hat N_1(\la) \cdots \hat N_{d'}(\la) = P^{(d')}(\la)\Lambda(\la),
$$
where for $i>0$, the matrices
$\hat N_i(\la):= (\Lambda_0 \cdots \Lambda_{i-1})^{-1} N_i (\Lambda_0 \cdots \Lambda_{i-1})$  
are unimodular, since
$$  (\Lambda_0 \cdots \Lambda_{i-1})^{-1} \left[\begin{array}{cc} I_{\rho_{i-1}} & Z_i \\ 0 & I_{\nu_{i-1}} \end{array}\right]
\left[\begin{array}{cc} I_{\rho_{i-1}} & 0 \\ 0 & Q_i \end{array}\right]
 (\Lambda_0 \cdots \Lambda_{i-1})$$
$$ =  (\Lambda_0 \cdots \Lambda_{i-1})^{-1} \left[\begin{array}{cc} I_{\rho_{i-1}} & Z_i \\ 0 & I_{\nu_{i-1}} \end{array}\right]   (\Lambda_0 \cdots \Lambda_{i-1}) \left[\begin{array}{cc} I_{\rho_{i-1}} &  0 \\ 0 & Q_i \end{array}\right] .
$$
 \hfill
\end{proof}

\subsection{Extracting compact decompositions} \label{Subsec:compact}

It is now easy to see that if we discard the last $n-r$ columns of the factorization \eqref{factor} and replace again $\la$ by $(\la-\la_0)$  then we obtain the form
\bigskip\begin{equation*} %\label{compactfactor}
P(\la) N_r(\la)= \hat M_r(\la)\Lambda_{(r)}(\la),  \quad \Lambda_{(r)}(\la)=\diag(I_{e_0}, (\la-\la_0) I_{e_1},\ldots,(\la-\la_0)^{d'} I_{e_{d'}}), 
\end{equation*}
where $N_r(\la)=N(\la)\begin{bsmallmatrix} I_{r} \\ 0 \end{bsmallmatrix}$ is a submatrix of a unimodular matrix,  $\hat M_r(\la)= P^{(d')}(\la)\begin{bsmallmatrix} I_r  \\ 0 \end{bsmallmatrix}$ is a polynomial matrix, 
and the constant matrix $\tilde M_r(\la_0)= L_0^{(d')}$ has full rank $r$.
This is the desired compact local Smith form
described in Theorem \ref{compactSM}. 

If we  apply the Toeplitz rank search algorithm to a rational matrix $R(\la) \in \C(\la)^{m\times n}$ without any poles at $\la_0$ and hence with a Taylor expansion
\begin{equation*} %\label{Taylor}  
R(\la) = R_0 +  R_{1}(\la-\la_0)^{1} +R_{2}(\la-\la_0)^{2} + R_{3}(\la-\la_0)^{3} + \ldots 
\end{equation*}
then we only need the leading terms $R_0, \ldots, R_{d'}$ of the expansion to obtain the same factorization as in \eqref{factor} except that the polynomial matrix $P^{(d')}(\la)$ is now replaced by a rational matrix $R^{(d')}(\la)$ whose Taylor expansion starts with the 
constant matrix $R^{(d')}(\la_0)=\left[L_0^{(d')}, \; 0 \;\right] $.
This leads to the following decomposition for a rational matrix $R(\la)$ which has zeros at $\la_0$ but no poles~:
\begin{equation*} %\label{compactrfactor}
R(\la) N_r(\la)= \hat M_r(\la)\Lambda_{(r)}(\la),  \quad \Lambda_{(r)}(\la)=\diag(I_{e_0}, (\la-\la_0) I_{e_1},\ldots,(\la-\la_0)^{d'} I_{e_{d'}}), 
\end{equation*}
where $N_r(\la)=N(\la)\begin{bsmallmatrix} I_{r} \\ 0 \end{bsmallmatrix}$ is a submatrix of a unimodular matrix,  $\hat M_r(\la)= R^{(d')}(\la)\begin{bsmallmatrix} I_r  \\ 0 \end{bsmallmatrix}$ is a rational matrix, and the constant matrix $\tilde M_r(\la_0)= L_0^{(d')}$ has full rank $r$.
This is the desired compact local Smith-McMillan form
described in Theorem \ref{compactMCM} for a rational matrix with a zero at $\la_0$ which is not a pole.

\medskip

We finally consider the case of a coalescent pole/zero. As pointed out in Theorem \ref{th:rank}, the Toeplitz matrices $T_{\la_0,k}(R)$
constructed with the coefficients of the Laurent expansion around the point $\la_0$
\begin{equation*} %\label{Laurent2}
R(\la) := R_{-\ell}(\la-\la_0)^{-\ell} + R_{-\ell +1}(\la-\la_0)^{-\ell+1} + R_{-\ell+2}
(\la-\la_0)^{-\ell+2} + R_{-\ell+3}(\la-\la_0)^{-\ell+3} + \ldots 
\end{equation*}
yields the structural indices $\sigma_i, i=1,\ldots,r$ of $R(\la)$ at the pole/zero $\la_0$. 
The way to reduce this to the rational case without a pole at $\la_0$ is to consider the scaled rational matrix $\widehat R(\la):= (\la-\la_0)^\ell R(\la)$. It is obvious that the structural indices $\widehat \sigma_i,\; i=1,\ldots,r$ of $\widehat R(\la)$ and $\sigma_i,\; i=1,\ldots,r$ of $R(\la)$ are related by a constant shift
$$   \widehat \sigma_i = \sigma_i + \ell \ge 0, \;\; i=1,\ldots,r.
$$
The Toeplitz rank search applied to $\widehat R(\la)$ then becomes a Taylor expansion of $\widehat R(\la)$, to which we can apply the results of the previous sections. 
After dividing $\widehat R(\la)$ and $\widehat \Lambda_{(r)}(\la)$ again by $(\la-\la_0)^\ell$, this leads to the following local decomposition for a general rational matrix $R(\la)$~:
\begin{equation*} %\label{compactrfactor2}
R(\la) N_r(\la)= \hat M_r(\la)\Lambda_{(r)}(\la),  \quad \Lambda_{(r)}(\la)=\diag((\la-\la_0)^{-\ell}I_{e_{-\ell}},\ldots,(\la-\la_0)^{d'} I_{e_{d'}}), 
\end{equation*}
where $N_r(\la)=N(\la)\begin{bsmallmatrix} I_{r} \\ 0 \end{bsmallmatrix}$ is a submatrix of a unimodular matrix,  $\hat M_r(\la)= R^{(d')}(\la)\begin{bsmallmatrix} I_r  \\ 0 \end{bsmallmatrix}$ is a rational matrix, and the constant matrix $\tilde M_r(\la_0)= L_0^{(d')}$ has full rank $r$.
This is the desired compact local Smith-McMillan form
described in Theorem \ref{compactMCM} for a general rational matrix with a pole/zero at $\la_0$.

\medskip

\section{Numerical examples} \label{Sec:numerics}

In this section we give a number of numerical results for the computation of the compact local Smith form
at the eigenvalue $\la_0=0$, computed using the algorithm\footnote{A MATLAB implementation of the algorithm we used is freely avaliable from github at the link \href{https://github.com/VanDoorenPaul/Compact-local-Smith-form}{https://github.com/VanDoorenPaul/Compact-local-Smith-form}}
described in Section \ref{Sec:Compact}. The test matrices were polynomial matrices of dimensions $4\times 5$ of normal rank 3 and with given invariant factors $\la^0, \la^1, \la^3$ at the eigenvalue $0$. The matrices were then constructed using the product
$$  P(\la)=M(\la)\left[ \begin{array}{ccc|cc} 1 & 0 & 0 & 0 & 0 \\ 0 & \la & 0 & 0 & 0 \\
0 & 0 & \la^3 & 0 & 0 \\ \hline 0 & 0 & 0 & 0 & 0 \\ 0 & 0 & 0 & 0 & 0 \end{array}\right] N(\la),
$$
where $M(\la)$ and $N(\la)$ are random polynomial matrices of respective dimensions 
$4\times 4$ and $5\times 5$, and of degree 2, which implies that the polynomial matrix $P(\la)$ has degree 8. The coefficients of the matrices $M(\la)$ and $N(\la)$ were generated using the $i$-th power of {\tt randn}, the random generator of Matlab with normal distribution. As a consequence, the dynamical range of the coefficients is growing and the norm of the matrix $P(\la)$ is typically growing as well with the power $i$. We used the Frobenius norm of a polynomial matrix $P(\la)$ of degree $d$, which is defined as follows~:
$$ \|P(\la)\|:=\|[ P_0, P_1, \ldots,P_d]\|_F.
$$
In Table \ref{table:T1} we show the results of our algorithm applied to the polynomial matrices $P(\la)$ generated for the powers $i$ going from 1 to 10. 

%\begin{center}
\begin{table}[h!]
    \centering
\begin{tabular}{rcccc}  
 $i$ &    $\|P\|$    & $\|ResP\|/\|P\|$ & $\|N\|$ \\ \hline
 1 &  2.5771e+01 &  2.5834e-16 &  4.2134e+00 \\
 2 &  5.3985e+01 &  4.8614e-16 &  6.0850e+00 \\
 3 &  4.0056e+02 &  1.0542e-15 &  5.9805e+01 \\
 4 &  1.8805e+03 &  5.2413e-15 &  7.6625e+01 \\
 5 &  5.6940e+03 &  2.9466e-15 &  3.6432e+01 \\
 6 &  1.2067e+04 &  2.0221e-16 &  3.9821e+00 \\
 7 &  2.4400e+04 &  2.4857e-15 &  2.7510e+01 \\
 8 &  2.1014e+05 &  2.8026e-11 &  3.7109e+03 \\
 9 &  2.5845e+05 &  1.6973e-13 &  5.1576e+03 \\
10 & 1.7714e+06 &  8.6854e-12 &  1.0296e+05  
\end{tabular}
 \caption{Recovery of the local Smith form
 for increasing powers $i$ of the random elements.}
    \label{table:T1}
\end{table}
%\end{center}

Column 2 and 4 give the Frobenius norms of the polynomial matrix $P(\la)$ and the unimodular matrix $N(\la)$. The accuracy of the computations is then verified using the Frobenius norm of the residual equation
$$ 
ResP(\la) = P(\la)N_r(\la)-\hat M_r(\la)\Lambda_{(r)}.
$$
Column 3 gives the relative norm $\|ResP\|/\|P\|$.
The structural indices were recovered correctly for all the test examples. 
It can be observed from these results that the accuracy of the algorithm is quite satisfactory, even for matrices with large dynamical range in the coefficients. The loss of accuracy, observed in some cases, is probably due the non-orthogonal Gram-Schmidt elimination step of our algorithm. But this could perhaps be improved by a single step of iterative refinement\cite{Abd71}.

In the second experiment, we check the robustness of our algorithm against polynomial matrices and Smith forms of high degree. We generated 10 matrices with local Smith form
$$  P(\la)=M(\la)\left[ \begin{array}{ccc|cc} 1 & 0 & 0 & 0 & 0 \\ 0 & \la^{k+1} & 0 & 0 & 0 \\
0 & 0 & \la^{k+2} & 0 & 0 \\ \hline 0 & 0 & 0 & 0 & 0 \\ 0 & 0 & 0 & 0 & 0 \end{array}\right] N(\la),
$$
for $k=1:10$, and transformation matrices $M(\la)$ and $N(\la)$ of degree 10. The degree of the polynomial matrices is therefore
$22+k$. The relative precision of the obtained decomposition
is again verified using the ratio $\|ResP\|/\|P\|$ of the Frobenius norm of the residual equation and the Frobenius norm of the matrix $P(\la)$. The structural indices were again recovered correctly for all the test examples.  This experiment shows that the method has remarkable stability properties, even for large degree polynomial matrices. 

%\begin{center} 
\begin{table}[h!]
\centering
\begin{tabular}{rcccc} 
 $k$ &    $\|P\|$    & $\|ResP\|/\|P\|$ & $\|N\|$ \\ \hline
1 &  7.1558e+01 &  2.7897e-16 &  3.0163e+00 \\
2 &  7.5268e+01 &  3.4156e-16 &  2.5936e+00 \\
3 &  7.6952e+01 &  2.0202e-16 &  1.8877e+00 \\
4 &  7.9385e+01 &  9.2658e-16 &  5.2263e+00 \\
5 &  7.9653e+01 &  2.6356e-16 &  2.1519e+00 \\
6 &  8.0100e+01 &  5.1580e-16 &  2.8901e+00 \\
7 &  8.1521e+01 &  6.5986e-16 &  4.5883e+00 \\
8 &  8.6119e+01 &  6.7279e-15 &  6.4173e+01 \\
9 &  8.7700e+01 &  5.8009e-16 &  3.7989e+00 \\
10 & 8.8879e+01 &  4.3347e-16 &  2.9266e+00
\end{tabular}
 \caption{Recovery of the local Smith form
 for increasing degrees of $P(\la)$.}
    \label{table:T2}
\end{table}
%\end{center}

\section{Conclusions} \label{Sec:conclusion}

In this paper we revisited the Toeplitz rank search algorithm developed in \cite{vdv} and showed that an appropriately modified variant also constructs a compact local Smith-McMillan decomposition at a given expansion point $\la_0$ that is a pole/zero of a rational matrix $R(\la)$. In this process we construct a unimodular transformation matrix 
whose columns are root polynomials introduced in \cite{GLR82} for regular polynomial matrices, refined in \cite{DopN21} for singular polynomial matrices, and extended in \cite{NofV22} to rational matrices. 
As a consequence, the degree of the constructed unimodular transformation matrix, is of minimal degree.
We also showed that, when applied to a polynomial matrix $P(\la)$, this decomposition is a compact local Smith form of $P(\la)$. The special type of transformation matrices used in this paper are not orthogonal, but are more related to a Gram-Schmidt orthogonalization procedure. This is reassuring since there exist numerically reliable implementations of the classical Gram-Schmidt procedure \cite{Abd71}, which might also apply to the Toeplitz rank search algorithm explained in this paper.  

\medskip

\appendix 

\renewcommand{\thesection}{A}

\section*{Appendix A}

\begin{theorem}\label{ohyes}
Fix an elementary divisor domain $\cal R$ and let $\F$ be the field of fractions of $\cal R$. Let $L\in {\cal R}^{r \times r}$ be invertible over $\F$, and let $X \in {\cal R}^{r \times (n-r)}$. Then the following are equivalent
\begin{enumerate}
\item  $\begin{bmatrix}
L&X
\end{bmatrix}$ is ${\cal R}$-unimodularly equivalent to $\begin{bmatrix}
L&0
\end{bmatrix}$
\item $\begin{bmatrix}
L&X\\
0&0
\end{bmatrix}$ is ${\cal R}$-unimodularly equivalent to $\begin{bmatrix}
L&0\\
0&0
\end{bmatrix}$;
\item $X=LY$ for some $Y \in {\cal R}^{r \times (n-r)}$.
\end{enumerate} 
\end{theorem}
\begin{proof}
\begin{itemize}
\item[$1 \Rightarrow 2$]There exist unimodular $U,V$ such that $U\begin{bmatrix}
L&X
\end{bmatrix}=\begin{bmatrix}
L&0
\end{bmatrix}V$, then
\[ (U \oplus I) \begin{bmatrix}
L&X\\
0&0
\end{bmatrix}=\begin{bmatrix}
L&0\\
0&0
\end{bmatrix}V.\] 
\item[$2 \Rightarrow 3$]  It holds
\[ \begin{bmatrix}
U_1 & U_2\\
U_3 & U_4
\end{bmatrix} \begin{bmatrix}
L&X\\
0&0
\end{bmatrix}=\begin{bmatrix}
L&0\\
0&0
\end{bmatrix}\begin{bmatrix}
V_1 & V_2\\
V_3 & V_4
\end{bmatrix} \Rightarrow \begin{cases}
U_3 L = 0 = U_3 X;\\
L V_1 = U_1 L;\\
L V_2 = U_1 X
\end{cases},\]
where the matrices whose blocks are $U_i$ and $V_i$ are both unimodular. But $U_3 L=0 \Rightarrow U_3=0$ because $L$ is invertible over $\F$. Thus $U_1$ must be unimodular, whence $X=U_1^{-1} L V_2 $. On the other hand $V_1 = L^{-1} U_1 L$ must be unimodular since $\det (V_1) = \det(L)^{-1} \det(U_1) \det(L) =\det(U_1)$ is a unit and $U_1^{-1} L = L V_1^{-1}$. Defining $Y=V_1^{-1} V_2$ we then have $X=LY$, as sought.
\item[$3 \Rightarrow 1$] We have
$$ \begin{bmatrix}
L & X
\end{bmatrix}=\begin{bmatrix}
L & 0
\end{bmatrix}\begin{bmatrix}
I & Y\\
0 & I
\end{bmatrix}.$$
\end{itemize}
\end{proof}

\begin{theorem} \label{ring}
Let ${\cal R}$ be an elementary divisor domain with field of fractions $\F$, and suppose that the Smith forms over ${\cal R}$ of $\begin{bmatrix}
N & C
\end{bmatrix}$ and $\begin{bmatrix}
N & 0
\end{bmatrix}$ are the same where $N \in {\cal R}^{m \times r}$ has full column rank and $C \in {\cal R}^{m \times (n-r)}$. Then, $C=NY$ for some $Y \in {\cal R}^{r \times (n-r)}$.
\end{theorem}
\begin{proof}
Since $N$ has full column rank we can write it as $N=U \begin{bmatrix}
T\\
0
\end{bmatrix}$ where $U \in {\cal R}^{m\times m}$ is unimodular and $T \in {\cal R}^{r \times r}$ (e.g. using the Hermite normal form \cite{F}). Moreover, $T$ is invertible over $\F$ because $\rank(T)=\rank(N)=r$. Hence,
$$ U^{-1} \begin{bmatrix}
N & C
\end{bmatrix} = \begin{bmatrix}
T & C_0\\
0 & C_1
\end{bmatrix};$$
but then $r \geq \rank(T) + \rank(C_1) \Rightarrow C_1=0$. By Lemma \ref{ohyes}, this implies in turn that $C_0=T Y$ for some $Y \in {\cal R}^{r \times (n-r)}$. However,
\[ \begin{bmatrix}
N & C
\end{bmatrix} = U \begin{bmatrix}
T & TY\\
0 & 0
\end{bmatrix} \Rightarrow C=NY.\]
\end{proof}

Lemma \ref{Lem:extend} then follows by applying Theorem \ref{ring} to the matrix $\begin{bmatrix}
    \pi(\la)R(\la) N_r(\la) & \pi(\la) R(\la)\tilde N(\la)
\end{bmatrix}.$

\end{document}